\documentclass[12pt,a4paper]{amsart}
\usepackage{amssymb,amscd}

\theoremstyle{plain}
\newtheorem{theorem}{Theorem}
\newtheorem{proposition}{Proposition}[section]
\newtheorem{lemma}[proposition]{Lemma}

\theoremstyle{definition}

\newtheorem{example}[proposition]{Example}

\numberwithin{equation}{section}

\newcommand\bA{{\mathbb A}}

\newcommand\bG{{\mathbb G}}

\newcommand\bP{{\mathbb P}}

\newcommand\cL{{\mathcal L}}
\newcommand\cM{{\mathcal M}}

\newcommand\cO{{\mathcal O}}

\newcommand\cT{{\mathcal T}}

\newcommand\fm{{\mathfrak m}}

\newcommand\aff{\operatorname{aff}}

\newcommand\id{\operatorname{id}}

\newcommand\op{\operatorname{op}}

\newcommand\Aut{\operatorname{Aut}}

\newcommand\GL{\operatorname{GL}}

\newcommand\Hom{\operatorname{Hom}}

\newcommand\Ima{\operatorname{Im}}
\newcommand\Lie{\operatorname{Lie}}

\newcommand\Pic{\operatorname{Pic}}

\newcommand\PGL{\operatorname{PGL}}

\newcommand\Supp{\operatorname{Supp}}

\title[Actions of non-affine algebraic groups]
{Some basic results on actions of non-affine algebraic groups}

\author{Michel Brion}
\address{Universit\'e de Grenoble I\\
D\'epartement de Math\'ematiques\\
Institut Fourier, UMR 5582 du CNRS\\
38402 Saint-Martin d'H\`eres Cedex, France}
\email{Michel.Brion@ujf-grenoble.fr}

\begin{document}
 
\begin{abstract}
We study actions of connected algebraic groups on normal algebraic
varieties, and show how to reduce them to actions of affine subgroups.
\end{abstract}

\maketitle

\section{Introduction}
\label{sec:introduction}

Algebraic group actions have been extensively studied under the
assumption that the acting group is affine or, equivalently, linear;
see \cite{KSS89,MFK94,PV94}. In contrast, little seems to be known
about actions of non-affine algebraic groups. In this paper, we show
that these actions may be reduced to actions of affine
subgroup schemes, in the setting of normal varieties. 

\smallskip

Our starting point is the following theorem of Nishi and Matsumura
(see \cite{Ma63}). Let $G$ be a connected algebraic group of
automorphisms of a nonsingular algebraic variety $X$ and denote by
$$
\alpha_X : X \longrightarrow A(X)
$$
the \emph{Albanese morphism}, that is, the universal morphism to an
abelian variety (see \cite{Se58b}). Then \emph{$G$ acts on $A(X)$ by
translations, compatibly with its action on $X$, and the kernel of the
induced homomorphism $G \to A(X)$ is affine}.

\smallskip

When applied to the action of $G$ on itself via left multiplication,
this shows that the Albanese morphism 
$$
\alpha_G : G \longrightarrow A(G)
$$ 
is a surjective group homomorphism having an affine kernel. Since this
kernel is easily seen to be smooth and connected, this gives back
Chevalley's structure theorem: \emph{any connected algebraic group $G$
is an extension of an abelian variety $A(G)$ by a connected affine
algebraic group} $G_{\aff}$ (see \cite{Co02}) for a modern proof). 

\smallskip

The Nishi--Matsumura theorem may be reformulated as follows:
for any faithful action of $G$ on a nonsingular variety $X$, the
induced homomorphism $G \to A(X)$ factors through a homomorphism 
$A(G) \to A(X)$ having a finite kernel (see \cite{Ma63} again). 
This easily implies the existence of a $G$-equivariant morphism
$$
\psi: X \longrightarrow A,
$$ 
where $A$ is an abelian variety, quotient of $A(G)$ by a finite
subgroup scheme (see Section \ref{sec:proof3} for details). 
Equivalently, $A \cong G/H$ where $H$ is a closed subgroup scheme of
$G$ such that $H \supset G_{\aff}$ and the quotient $H/G_{\aff}$ is
finite; in particular, $H$ is affine, normalized by $G$, and uniquely
determined by $A$. Then there is a $G$-equivariant isomorphism
$$
X \cong G \times^H Y,
$$ 
where the right-hand side denotes the homogeneous fiber bundle over
$G/H$ associated to the scheme-theoretic fiber $Y$ of $\psi$ at the
base point.

\smallskip

In particular, given a faithful action of an abelian variety $A$ on a
nonsingular variety $X$, there exist a positive integer $n$ and a
closed $A_n$-stable subscheme $Y \subset X$  such that 
$X \cong A \times^{A_n} Y$, where $A_n\subset A$ denotes the kernel of
the multiplication by $n$. For free actions (that is, abelian
torsors), this result is due to Serre, see \cite[Prop.~17]{Se58a}.

\smallskip

Next, consider a faithful action of $G$ on a possibly singular variety
$X$. Then, besides the Albanese morphism, we have the \emph{Albanese map}
$$
\alpha_{X,r}: X - - \to A(X)_r,
$$
i.e., the universal rational map to an abelian variety. Moreover,
the regular locus $U \subset X$ is $G$-stable, and
$A(U) = A(U)_r = A(X)_r$. Thus, $G$ acts on $A(X)_r$ via a
homomorphism $A(G) \to A(X)_r$ such that the canonical homomorphism 
$$
h_X : A(X)_r \longrightarrow A(X)
$$ 
is equivariant; $h_X$ is surjective, but generally not an isomorphism,
see \cite{Se58b} again. Applying the Nishi--Matsumura theorem to $U$,
we see that the kernel of the $G$-action on $A(X)_r$ is affine, and
there exists a $G$-equivariant rational map
$\psi_r : X - - \to A$
for some abelian variety $A$ as above. 

\smallskip

However, $G$ may well act trivially on $A(X)$; then there exists no
morphism $\psi$ as above, and $X$ admits no equivariant embedding into
a nonsingular $G$-variety. This happens for several classes of
examples constructed by Raynaud, see \cite[XII 1.2, XIII 3.2]{Ra70} or
Examples \ref{ex:raynaud1}, \ref{ex:raynaud2}, \ref{ex:raynaud3}; in
the latter example, $X$ is normal, and $G$ is an abelian variety
acting freely. However, we shall show that such an equivariant embedding
(in particular, such a morphism $\psi$) exists locally for any normal
$G$-variety.

\smallskip

To state our results in a precise way, we introduce some notation and 
conventions. We consider algebraic varieties and schemes over an
algebraically closed field $k$; morphisms are understood to be
$k$-morphisms. By a variety, we mean a separated integral scheme of
finite type over $k$; a point will always mean a closed point.
As a general reference for algebraic geometry, we use the book
\cite{Ha77}, and \cite{DG70} for algebraic groups.

\smallskip

We fix a connected algebraic group $G$, that is, a $k$-group
variety; in particular, $G$ is smooth. A $G$-\emph{variety}
is a variety $X$ equipped with an algebraic $G$-action
$$
\varphi : G \times X \longrightarrow X, \quad (g,x) \longmapsto g \cdot x.
$$ 
The \emph{kernel} of $\varphi$ is the largest subgroup scheme of $G$
that acts trivially. We say that $\varphi$ is \emph{faithful} if its 
kernel is trivial. The $G$-variety $X$ is \emph{homogeneous} 
(resp.~\emph{almost homogeneous}) if it contains a unique orbit
(resp.~an open orbit). Finally, a $G$-\emph{morphism} is an
equivariant morphism between $G$-varieties.

\smallskip

We may now formulate our results, and discuss their relations to
earlier work.

\begin{theorem}\label{thm:quasi-projective}
Any normal $G$-variety admits an open covering by $G$-stable
quasi-projective subsets.
\end{theorem}

When $G$ is affine, this fundamental result is due to Sumihiro (see
\cite{KSS89, Su74}); on the other hand, it has been obtained by
Raynaud in another setting, namely, for actions of smooth, connected
group schemes on smooth schemes over normal bases (in particular, for
actions of connected algebraic groups on nonsingular varieties), see 
\cite[Cor.~V 3.14]{Ra70}.

\smallskip

Theorem \ref{thm:quasi-projective} implies readily the
quasi-projectivity of homogeneous varieties and, more generally, of
normal varieties where a connected algebraic group acts with a unique
closed orbit. For the latter result, the normality assumption cannot
be omitted, as shown by Example \ref{ex:raynaud1}.

\smallskip

Next, we obtain a version of the Nishi--Matsumura theorem: 

\begin{theorem}\label{thm:fibration}
Let $X$ be a normal, quasi-projective variety on which $G$ acts
faithfully. Then there exists a $G$-morphism
$$
\psi : X \longrightarrow A,
$$ 
where $A$ is the quotient of $A(G)$ by a finite subgroup scheme.
Moreover, $X$ admits a $G$-equivariant embedding into the
projectivization of a $G$-homogeneous vector bundle over $A$.
\end{theorem}

The second assertion generalizes (and builds on) another result of
Sumihiro: when $G$ is affine, any normal quasi-projective $G$-variety
admits an equivariant embedding into the projectivization of a
$G$-module (see \cite{KSS89, Su74} again). It implies that any normal,
quasi-projective $G$-variety has an equivariant embedding into a
nonsingular $G$-variety, of a very special type. Namely, a vector
bundle $E$ over an abelian variety $A$ is homogeneous with
respect to a connected algebraic group $G$ (acting transitively on $A$)
if and only if $a^*(E) \cong E$ for all $a \in A$; see \cite{Mu78} for
a description of all homogeneous vector bundles on an abelian variety.
 
\smallskip

Clearly, a morphism $\psi$ as in Theorem \ref{thm:fibration} is
not unique, as we may replace $A$ with any finite quotient. In
characteristic zero, we can always impose that the fibers of $\psi$
are normal varieties, by using the Stein factorization. Under that
assumption, there may still exist several morphisms $\psi$, but no
universal morphism (see Example \ref{ex:several}). 

\smallskip

Returning to arbitary characteristics, there does exist a universal
morphism $\psi$ when $X$ is normal and almost homogeneous, namely, the
Albanese morphism:

\begin{theorem}\label{thm:albanese}
Let $X$ be a normal, almost homogeneous $G$-variety. Then 
$A(X) = A(X)_r = G/H$, where $H \subset G$ is a closed subgroup scheme
containing $G_{\aff}$, and the quotient group scheme $H/G_{\aff}$ is
finite. Moreover, each fiber of the Albanese morphism is a normal
variety, stable under $H$ and almost homogeneous under $G_{\aff}$.
\end{theorem}

This \emph{Albanese fibration} is well-known in the setting of complex
Lie groups acting on compact K\"ahler manifolds (the Remmert--van de
Ven theorem, see e.g. \cite[Sec.~3.9]{Ak95}), and easily obtained for
nonsingular varieties (see \cite[Sec.~2.4]{Br06}). 

\smallskip

In particular, Theorem \ref{thm:albanese} applies to any normal
\emph{$G$-equivariant embedding} of $G$, i.e., to a normal $G$-variety
$X$ containing an open orbit isomorphic to $G$: then 
$A(X) = A(G) = G/G_{\aff}$ and hence
$$
X \cong G \times^{G_{\aff}} Y,
$$
where $Y$ is a normal $G_{\aff}$-equivariant embedding of $G_{\aff}$.
(Here again, the normality assumption cannot be omitted, see Example
\ref{ex:raynaud1}). If $G$ is a \emph{semi-abelian variety}, that is,
$G_{\aff}$ is a torus, then $Y$ is a toric variety, and the
$G$-equivariant embedding $X$ is called a \emph{semi-abelic variety}. 
In that case, the above isomorphism has been obtained by Alexeev under
the assumption that $X$ is projective (see \cite[Sec.~5]{Al02}). As
another application, consider a normal algebraic monoid $X$ with unit
group $G$; then $X$ may be regarded as a $(G \times G$)-equivariant
embedding of $G$, and the above isomorphism yields another proof of
the main result of \cite{BR06}.

\smallskip

For a complete homogeneous variety, it is known that the Albanese
fibration is trivial. In the setting of compact homogeneous K\"ahler
manifolds, this is the classical Borel--Remmert theorem (see
e.g. \cite[Sec.~3.9]{Ak95}); in our algebraic setting, this is part of
a structure result due to Sancho de Salas:

\begin{theorem}\label{thm:product}
{\rm (\cite[Thm.~5.2]{Sa03})}. Let $X$ be a complete variety,
homogeneous under a faithful action of $G$. Then there exists a
canonical isomorphism $G \cong A(G) \times G_{\aff}$, and $G_{\aff}$
is semi-simple of adjoint type. Moreover, there exists a canonical
isomoporhism of $G$-varieties $X \cong A(G) \times Y$, where $Y$ is a
complete, homogeneous $G_{\aff}$-variety; the Albanese morphism of $X$
is the first projection $X \to A(G)$. 
\end{theorem}

\smallskip

Here we present a short proof of that result by analyzing first the
structure of $G$ and then the Albanese morphism of $X$, while
\cite{Sa03} proceeds by constructing the second projection  
$X \to Y$ via the study of the $G_{\aff}$-action on $X$. In
characteristic zero, this projection is the \emph{Tits morphism} that
assigns to any point of $X$ its isotropy Lie algebra; this yields a
very simple proof of Theorem \ref{thm:product}, see
\cite[Sec.~1.4]{Br06}. But this approach does not extend to positive
characteristics, as the $G$-action on $X$ may well be non-smooth;
equivalently, the orbit maps may be non-separable. In fact, there
exist many complete varieties $Y$ that are homogeneous under a
non-smooth action of a semi-simple group $G_{\aff}$. These 
\emph{varieties of unseparated flags} are classified in 
\cite{HL93, Sa03, We93}. Any such variety $Y$ is projective (e.g., by
Theorem \ref{thm:quasi-projective}) and rational. Indeed, $Y$ contains
only finitely many fixed points of a maximal torus of $G_{\aff}$;
thus, $Y$ is paved by affine spaces, see \cite{Bi73}.

\smallskip

In positive characteristics again, the homogeneity assumption of
Theorem \ref{thm:product} cannot be replaced with the assumption that
the tangent bundle is globally generated. Indeed, there exist
nonsingular projective varieties which have a trivial tangent bundle,
but fail to be homogeneous (see \cite{Ig55} or Example
\ref{ex:igusa}). Also, for $G$-varieties of unseparated flags,
the subbundle of the tangent bundle generated by $\Lie(G)$ is a proper
direct summand (see Example \ref{ex:split}). 

\smallskip

The proofs of Theorems 1--4 are given in the corresponding
sections. The final Section \ref{sec:examples} presents examples
illustrating the assumptions of these theorems.  Three of these
examples are based on constructions, due to Raynaud, of torsors under
abelian schemes, that we have adapted so as to make them known to
experts in algebraic groups.

\smallskip

In the opposite direction, it would be interesting to extend our
results to group schemes. Our methods yield insight into actions of
abelian schemes (e.g., the proof of Theorem \ref{thm:fibration} adapts 
readily to that setting), but the general case requires new ideas.

\medskip

\noindent
{\bf Acknowledgements.} I wish to thank Jos\'e Bertin, St\'ephane
Druel, Adrien Dubouloz, and Ga\"el R\'emond for many useful
discussions, and Patrick Brosnan for pointing out the reference
\cite{Sa03}.

\section{Proof of Theorem \ref{thm:quasi-projective}}
\label{sec:proof1}

As in the proof of the projectivity of abelian varieties (see
e.g. \cite[Thm.~7.1]{Mi86}), we first obtain a version of the theorem
of the square. We briefly present the setting of this theorem,
refererring to \cite[Chap.~IV]{Ra70} and \cite[Sec.~6.3]{BLR90} for
further developments. 

Consider an invertible sheaf $L$ on a $G$-variety $X$. For any 
$g \in G$, let $g^*(L)$ denote the pull-back of $L$ under the
automorphism of $X$ induced by $g$. In loose words, $L$ satisfies the
theorem of the square if there are compatible isomorphisms
\begin{equation}\label{eqn:square}
(gh)^*(L) \otimes L \cong g^*(L) \otimes h^*(L)
\end{equation}
for all $g,h$ in $G$. Specifically, consider the action
$\varphi : G \times X \to X$, the second projection 
$p_2: G \times X \to X$, and put 
\begin{equation}\label{eqn:cL}
\cL := \varphi^*(L) \otimes p_2^*(L)^{-1}.
\end{equation}
This is an invertible sheaf on $G \times X$, satisfying
$$
\cL \vert_{\{g\} \times X} \cong g^*(L) \otimes L^{-1}
$$
for all $g\in G$. Next, consider the variety $G \times G \times X$ over
$X$, and the $X$-morphisms
$$
m, p_1,p_2 : G \times G \times X \longrightarrow G \times X 
$$
given by the multiplication and the projections 
$G \times G \to G$. Then 
$$
\cM := m^*(\cL) \otimes p_1^*(\cL)^{-1} \otimes p_2^*(\cL)^{-1}
$$
is an invertible sheaf on $G \times G \times X$, such that
$$
\cM \vert_{\{ (g,h) \} \times X} \cong (gh)^*(L) \otimes L \otimes
g^*(L)^{-1} \otimes h^*(L)^{-1}
$$ 
for all $g, h$ in $G$. Now $L$ satisfies the theorem of the square, if
$\cM$ is the pull-back of some invertible sheaf under the projection 
$$
f : G \times G \times X \longrightarrow G \times G.
$$
Then each $\cM \vert_{\{ (g,h) \} \times X}$ is trivial, which implies
(\ref{eqn:square}).

Also, recall the classical notion of a $G$-\emph{linearization} of the
invertible sheaf $L$, that is, a $G$-action of the total space of the
associated line bundle which is compatible with the $G$-action on
$X$ and commutes with the natural $\bG_m$-action (see 
\cite{KSS89, MFK94}). The isomorphism classes of
$G$-linearized invertible sheaves form a group denoted by
$\Pic^G(X)$. We will use the following observation (see
\cite[p.~32]{MFK94}):
 
\begin{lemma}\label{lem:descent}
Let $\pi: X \to Y$ be a torsor under $G$ for the fppf topology. Then
the pull-back under $\pi$ yields an isomorphism 
$\Pic(Y) \cong \Pic^G(X)$.  
\end{lemma}

\begin{proof}
By assumption, $\pi$ is faithfully flat and fits into a cartesian
square
$$
\CD
G \times X @>{\varphi}>> X \\
@V{p_2}VV @V{\pi}VV \\
X @>{\pi}>> Y. \\
\endCD
$$
Moreover, a $G$-linearization of an invertible sheaf $L$ on $X$ is
exactly a descent datum for $L$ under $\pi$, see
\cite[Sec.~1.3]{MFK94}. So the assertion follows from faithfully flat
descent, see e.g. \cite[Sec.~6.1]{BLR90}. 
\end{proof}

We now come to our version of the theorem of the square:

\begin{lemma}\label{lem:square}
Let $L$ be a $G_{\aff}$-linearizable invertible sheaf on a $G$-variety
$X$. Then $L$ satisfies the theorem of the square.
\end{lemma}

\begin{proof}
Consider the action of $G$ on $G \times X$ via left multiplication on
the first factor. Then $\varphi$ is equivariant, and $p_2$ is
invariant. Hence the choice of a $G_{\aff}$-linearization of $L$
yields a $G_{\aff}$-linearization of the invertible sheaf $\cL$ on 
$G \times X$ defined by (\ref{eqn:cL}). Note that the map
$$
\alpha_G \times \id_X : G \times X \longrightarrow  A(G) \times X
$$
is a $G_{\aff}$-torsor. By Lemma \ref{lem:descent}, there exists a
unique invertible sheaf $\cL'$ on $A(G) \times X$ such that
\begin{equation}\label{eqn:cM}
\cL = (\alpha_G \times \id_X)^*(\cL')
\end{equation}
as $G_{\aff}$-linearized sheaves. Then
$\cM = (\alpha_G \times \alpha_G \times \id_X)^*(\cM')$, where
$$
\cM' := m'{}^*(\cL') \otimes p'_1{}^*(\cL')^{-1} \otimes p'_2{}^*(\cL')^{-1}
$$
and $m',p'_1,p'_2: A(G) \times A(G) \times X \to A(G) \times X$ are
defined analogously to $m,p_1,p_2$. Thus, it suffices to show that
$\cM'$ is the pull-back of an invertible sheaf under the projection 
$$
f': A(G) \times A(G) \times X \longrightarrow A(G) \times A(G).
$$

Choose $x \in X$ and put 
$$
\cM'_x := \cM' \vert_{A(G) \times A(G) \times \{ x \}}.
$$
We consider $\cM'_x$ as an invertible sheaf on $A(G) \times A(G)$, and 
show that the invertible sheaf $\cM' \otimes f'{}^*(\cM'_x)^{-1}$ is
trivial. By a classical rigidity result (see \cite[Thm.~6.1]{Mi86}),
it suffices to check the triviality of the restrictions of this sheaf
to $\{ 0 \}\times A(G) \times X$, $A(G) \times \{ 0 \} \times X$, and 
$A(G) \times A(G) \times \{ x \}$. In view of the definition of
$\cM'_x$, it suffices in turn to show that
\begin{equation}\label{eqn:trivial}
\cM' \vert_{\{ 0 \} \times A(G) \times X} \cong \cO_{A(G) \times X}.
\end{equation}
For this, note that 
$\cL \vert_{G_{\aff} \times X} \cong \cO_{G_{\aff} \times X}$
since $L$ is $G_{\aff}$-linearized. By Lemma \ref{lem:descent}, it
follows that
\begin{equation}\label{eqn:rigidification}
\cL' \vert_{\{ 0 \} \times X} \cong \cO_{\{ 0 \} \times X}.
\end{equation}
Thus, $p'_1{}^*(\cL')$ is trivial; on the other hand,
$m' = p'_2$ on $\{ 0 \} \times A(G) \times X$. These facts imply 
(\ref{eqn:trivial}).
\end{proof}

Next, recall that for any invertible sheaf $L$ on a normal $G$-variety
$X$, some positive power $L^n$ admits a $G_{\aff}$-linearization; 
specifically, the Picard group of $G_{\aff}$ is finite, and we may
take for $n$ the order of that group (see \cite[p.~67]{KSS89}). 
Together with Lemma \ref{lem:square}, this yields: 

\begin{lemma}\label{lem:power}
Let $L$ be an invertible sheaf on a normal $G$-variety. Then some
positive power $L^n$ satisfies the theorem of the square; we may take
for $n$ the order of $\Pic(G_{\aff})$.
\end{lemma}

From this, we deduce an ampleness criterion on normal $G$-varieties,
analogous to a result of Raynaud about actions of group schemes on
smooth schemes over a normal base (see \cite[Thm.~V 3.10]{Ra70}):
 
\begin{lemma}\label{lem:divisor}
Let $X$ be a normal $G$-variety, and $D$ an effective Weil divisor on
$X$. If $\Supp(D)$ contains no $G$-orbit, then some positive multiple
of $D$ is a Cartier divisor generated by its global sections. If, in
addition, the open subset $X \setminus \Supp(D) \subset X$ is affine,
then some positive multiple of $D$ is ample.
\end{lemma}

\begin{proof}
Consider the regular locus $X_0 \subset X$ and the restricted divisor
$D_0 := D \cap X_0$. Then $X_0$ is $G$-stable, and the sheaf
$\cO_{X_0}(D_0)$ is invertible; moreover,
$g^*(\cO_{X_0}(D_0)) = \cO_{X_0}(g \cdot D_0)$
for all $g \in G$. By Lemma \ref{lem:power}, there exists a
positive integer $n$ such that $\cO_{X_0}(nD_0)$ satisfies the theorem
of the square. Replacing $D$ with $nD$, we obtain isomorphisms
$\cO_{X_0}(2D_0) \cong \cO_{X_0}(g \cdot D_0 + g^{-1} \cdot  D_0)$
for all $g \in G$. Since $X$ is normal, it follows that 
\begin{equation}\label{eqn:2D}
\cO_X(2D) \cong \cO_X(g \cdot D + g^{-1} \cdot D)
\end{equation}
for all $g \in G$. 

Let $U := X \setminus \Supp(D)$. By (\ref{eqn:2D}), the Weil divisor
$2D$ restricts to a Cartier divisor on every open subset
$$
V_g := X \setminus \Supp(g \cdot D + g^{-1} \cdot D) = 
g \cdot U \cap g^{-1} \cdot U, \quad g \in G.
$$
Moreover, these subsets form a covering of $X$ (indeed, given 
$x \in X$, the subset
$$
W_x := \{ g\in G ~\vert~ g \cdot x \in U\} \subset G
$$ 
is open, and non-empty since $U$ contains no $G$-orbit. Thus, $W_x$
meets its image under the inverse map of $G$, that is, there exists 
$g \in G$ such that $U$ contains both points $g \cdot x$ and 
$g^{-1} \cdot x$). It follows that $2D$ is a Cartier divisor on $X$.
Likewise, the global sections of $\cO_X(2D)$ generate this divisor on 
each $V_g$, and hence everywhere.

If $U$ is affine, then each $V_g$ is affine as well. Hence $X$ is
covered by affine open subsets $X_s$, where $s$ is a global section
of $\cO_X(2D)$. Thus, $2D$ is ample.
\end{proof}

We may now prove Theorem \ref{thm:quasi-projective}. 
Let $x \in X$ and choose an affine open subset $U$ containing
$x$. Then $G \cdot U$ is a $G$-stable open subset of $X$; moreover,
the complement $(G \cdot U) \setminus U$ is of pure codimension $1$
and contains no $G$-orbit. Thus, $G \cdot U$ is quasi-projective by
Lemma \ref{lem:divisor}.

\section{Proof of Theorem \ref{thm:fibration}}
\label{sec:proof2}

First, we gather preliminary results about algebraic groups and
their actions.

\begin{lemma}\label{lem:affine}
{\rm (i)} Let $\pi: X \to Y$ be a torsor under a group scheme
$H$. Then $\pi$ is affine if and only if $H$ is affine.

\smallskip

\noindent
{\rm (ii)} 
Let $X$ be a variety on which $G$ acts faithfully. Then the orbit map 
$$
\varphi_x : G \longrightarrow G \cdot x, \quad g \longmapsto g \cdot x
$$
is affine, for any $x \in X$.

\smallskip

\noindent
{\rm (iii)} Let $C(G)$ denote the center of $G$, and $C(G)^0$ its
reduced neutral component. Then $G = C(G)^0 G_{\aff}$.

\smallskip

\noindent
{\rm (iv)} Let $M$ be an invertible sheaf on $A(G)$, and 
$L := \alpha_G^*(M)$ the corresponding $G_{\aff}$-linearized
invertible sheaf on $G$. Then $L$ is ample if and only if $M$ is
ample. 
\end{lemma}

\begin{proof}
(i) follows by faithfully flat descent, like Lemma \ref{lem:descent}.

(ii) Since $\varphi_x$ is a torsor under the isotropy subgroup scheme
$G_x \subset G$, it suffices to check that $G_x$ is affine or, 
equivalently, admits an injective representation in a 
finite-dimensional $k$-vector space. Such a representation is afforded
by the natural action of $G_x$ on the quotient $\cO_{X,x}/\fm_x^n$ for
$n \gg 0$, where $\cO_{X,x}$ denotes the local ring of $X$ at $x$,
with maximal ideal $\fm_x$; see \cite[Lem.~p.~154]{Ma63} for details.

(iii) By \cite[Lem.~2]{Se58a}, $\alpha_G$ restricts to a surjective
morphism $C \to A(G)$. Hence the restriction 
$C^0 \to A(G) = G/G_{\aff}$ is surjective as well.

(iv) Since the morphism $\alpha_G$ is affine, the ampleness of $M$
implies that of $L$. The converse holds by 
\cite[Lem.~XI 1.11.1]{Ra70}.
\end{proof}

Next, we obtain our main technical result:

\begin{lemma}\label{lem:equivalent}
Let $X$ be a variety on which $G$ acts faithfully. Then the
following conditions are equivalent: 

\smallskip

\noindent
{\rm (i)} There exists a $G$-morphism $\psi: X \to A$, where $A$ is
the quotient of $A(G)$ by a finite subgroup scheme.

\smallskip

\noindent
{\rm (ii)} There exists a $G_{\aff}$-linearized invertible
sheaf $L$ on $X$ such that $\varphi_x^*(L)$ is ample for any 
$x \in X$.

\smallskip

\noindent
{\rm (iii)} There exists a $G_{\aff}$-linearized invertible
sheaf $L$ on $X$ such that $\varphi_{x_0}^*(L)$ is ample for some
$x_0 \in X$.
\end{lemma}

\begin{proof}
(i) $\Rightarrow$ (ii) Choose an ample invertible sheaf $M$ on the
abelian variety $A$. We check that $L := \psi^*(M)$ satisfies
the assertion of (ii). Indeed, by the universal property of the
Albanese morphism $\alpha_G$, there exists a unique $G$-morphism
$\alpha_x : A(G) \to A$ such that the square
$$
\CD
G @>{\varphi_x}>> X   \\
@V{\alpha_G}VV @V{\psi}VV    \\
A(G) @>{\alpha_x}>> A \\
\endCD
$$
is commutative. By rigidity (see e.g. \cite[Cor.~3.6]{Mi86}),
$\alpha_x$ is the composite of a homomorphism and a translation. 
On the other hand, $\alpha_x$ is $G$-equivariant; thus, it is the
composite of the quotient map $A(G) \to A$ and a translation.
So $\alpha_x$ is finite, and hence $\alpha_x^*(M)$ is ample. By Lemma
\ref{lem:affine} (iv), it follows that
$\varphi_x^*(L) =\alpha_G^*(\alpha_x^*(M))$ is ample.

(ii) $\Rightarrow$ (iii) is left to the reader.

(iii) $\Rightarrow$ (i) Consider the invertible sheaves $\cL$ on 
$G \times X$, and $\cL'$ on $A(G) \times X$, defined by (\ref{eqn:cL})
and (\ref{eqn:cM}). Recall from (\ref{eqn:rigidification}) that $\cL'$
is equipped with a rigidification along $\{ 0 \} \times X$. By
\cite[Sec.~8.1]{BLR90}, it follows that $\cL'$ defines a morphism
$$
\psi: X \longrightarrow \Pic A(G), \quad 
x \longmapsto \cL'\vert_{A(G) \times \{x \}},
$$
where $\Pic A(G)$ is equipped with its scheme structure (reduced,
locally of finite type) for which the connected components are exactly
the cosets of the dual abelian variety, $A(G)^{\vee}$. Since $X$ is
connected, its image under $\psi$ is contained in a unique coset.

By construction, $\psi$ maps every point $x \in X$ to the isomorphism
class of the unique invertible sheaf $M_x$ on $A(G)$ such that  
$$
\alpha_G^*(M_x) = \varphi_x^*(L)
$$
as $G_{\aff}$-linearized sheaves on $G$. When $x = x_0$, the
invertible sheaf $\varphi_x^*(L)$ is ample, and hence
$M_x$ is ample by Lemma \ref{lem:affine} (iv). It follows that the
points of the image of $\psi$ are exactly the ample classes
$a^*(M_{x_0})$, where $a \in A(G)$.

Moreover, $\varphi_{g \cdot x} = \varphi_x \circ \rho(g)$
for all $g \in G$ and $x \in X$, where $\rho$ denotes right
multiplication in $G$. Thus,
$\alpha_G^*(M_{g\cdot x}) = \rho(g)^*( \alpha_G^*(M_x))$, 
that is, $\psi(g \cdot x) = \alpha_G(g)^* \psi(x)$. 
In other words, $\psi$ is $G$-equivariant, where $G$ acts on 
$\Pic A(G)$ via the homomorphism $\alpha_G : G \to A(G)$ and the
$A(G)$-action on $\Pic A(G)$ via pull-back. Hence the image of
$\psi$ is the $G$-orbit of $\psi(x_0)$, that is, the quotient $A(G)/F$
where $F$ denotes the scheme-theoretic kernel of the polarization
homomorphism
$$
A(G) \longrightarrow A(G)^{\vee}, \quad 
a \longmapsto a^*(M_{x_0}) \otimes M_{x_0}^{-1}.
$$
\end{proof}

We are now in a position to prove Theorem \ref{thm:fibration}.
Under the assumptions of that theorem, we may choose an ample
invertible sheaf $L$ on $X$. Then $\varphi_x^*(L)$ is ample for any 
$x \in X$, as follows from Lemma \ref{lem:affine} (ii). Moreover,
replacing $L$ with some positive power, we may assume that $L$ is
$G_{\aff}$-linearizable. Now Lemma \ref{lem:equivalent} yields a
$G$-morphism
$$
\psi: X \longrightarrow A,
$$
where $A \cong G/H$ for an affine subgroup scheme $H \subset G$. 
So $X \cong G \times^H Y$, where $Y\subset X$ is a closed $H$-stable
subscheme. To complete the proof, it suffices to show that $X$ (and
hence $Y$) admits an $H$-equivariant embedding into the
projectivization of an $H$-module. Equivalently, $X$ admits an ample,
$H$-linearized invertible sheaf. We already know this when $H$ is
smooth and connected; the general case may be reduced to that one, as
follows. 

We may consider $H$ as a closed subgroup scheme of some $\GL_n$.
Then $G_{\aff}$ embeds into $\GL_n$ as the reduced neutral component
of $H$. The homogeneous fiber bundle 
$\GL_n \times ^{G_{\aff}} X$ is a normal $\GL_n$-variety, see
\cite[Prop.~4]{Se58a}; it is quasi-projective by
\cite[Prop.~7.1]{MFK94}. The finite group scheme $H/G_{\aff}$ acts on
that variety, and this action commutes with the $\GL_n$-action. Hence
the quotient
$$
(\GL_n \times ^{G_{\aff}} X)/(H/G_{\aff}) = \GL_n \times ^H X
$$
is a normal, quasi-projective variety as well (see e.g.
\cite[\S~12]{Mu70}). So we my choose an ample,
$\GL_n$-linearized invertible sheaf $L$ on that quotient. The
pull-back of $L$ to $X \subset \GL_n \times ^H X$ is the desired
ample, $H$-linearized invertible sheaf.

\section{Proof of Theorem \ref{thm:albanese}}
\label{sec:proof3}

We begin with some observations about the Albanese morphism of a
$G$-variety $X$, based on the results of \cite{Se58b}. Since the
Albanese morphism of $G \times X$ is the product morphism 
$\alpha_G \times \alpha_X$, there exists a unique morphism of
varieties $\beta: A(G) \times A(X) \to A(X)$
such that the following square is commutative: 
$$
\CD
G \times X @>{\varphi}>> X \\
@V{\alpha_G \times \alpha_X}VV @V{\alpha_X}VV \\
A(G) \times A(X) @>{\beta}>> A(X). \\ 
\endCD
$$
Then $\beta$ is the composite of a homomorphism and a
translation. Moreover, $\beta(0,z) = z$ for all $z$ in the image of
$\alpha_X$. Since this image is not contained in a translate of a
smaller abelian variety, it follows that $\beta(0,z) = z$ for all 
$z \in A(X)$. Thus, $\beta(y,z) = \alpha_{\varphi}(y) + z$, where 
$$
\alpha_{\varphi}: A(G) \longrightarrow A (X)
$$ 
is a homomorphism. In other words, the $G$-action on $X$ induces an
action of $A(G)$ on $A(X)$ by translations via $\alpha_{\varphi}$.

Given an abelian variety $A$, the datum of a morphism $X \to A$ up to
a translation in $A$ is equivalent to that of a homomorphism 
$A(X) \to A$. Thus, the datum of a $G$-equivariant morphism  
$$
\psi : X \longrightarrow A
$$ 
up to a translation in $A$, is equivalent to that of a homomorphism
$$
\alpha_{\psi} : A(X) \longrightarrow A = A(G)/F
$$
such that the composite 
$\alpha_{\psi} \circ \alpha_{\varphi}$ equals the quotient morphism 
$$
q : A(G) \to A(G)/F
$$ 
up to a translation. Then the kernel of $\alpha_{\varphi}$ is
contained in $F$, and hence is finite.

Conversely, if $\alpha_{\varphi}$ has a finite kernel, then there exist
a finite subgroup scheme $F \subset A(G)$ and a homomorphism
$\alpha_{\psi}$ such that $\alpha_{\psi} \circ \alpha_{\varphi} = q$.
Indeed, this follows from Lemma \ref{lem:equivalent} applied to the
image of $\alpha_{\varphi}$; alternatively, this is a consequence of
the Poincar\'e complete reducibility theorem (see e.g.  
\cite[Prop.~12.1]{Mi86}).
We also see that for a fixed subgroup scheme $F \subset A(G)$ (or,
equivalently, $H \subset G$), the set of homomorphisms $\alpha_{\psi}$ 
is a torsor under $\Hom(A(X)/\alpha_{\varphi}(A(G)), A(G)/F)$, a free
abelian group of finite rank.

So we have shown that the existence of a $G$-morphism $\psi$ as in
Theorem \ref{thm:fibration} is equivalent to the kernel of the
$G$-action on $A(X)$ being affine, and also to the kernel of
$\alpha_{\varphi}$ being finite.

Next, we assume that $X$ is normal and quasi-homogeneous, and prove
Theorem \ref{thm:albanese}. Choose a point $x$ in the open orbit
$X_0$ and denote by $G_x \subset G$ its isotropy subgroup scheme, so 
that $X_0$ is isomorphic to the homogeneous space $G/G_x$. Then $G_x$
is affine by Lemma \ref{lem:affine}; therefore, the product 
$H_0 := G_x G_{\aff}$ is a closed affine subgroup scheme of $G$, and
the quotient $F := H_0/G_{\aff}$ is finite. As a consequence, the
homogeneous space
$$
G/H_0 \cong (G/G_{\aff}) /(H_0/G_{\aff}) = A(G)/F
$$
is an abelian variety, isogenous to $A(G)$.

We claim that the Albanese morphism of $X_0$ is the quotient map 
$G/G_x \to G/H_0$. Indeed, let $f: X_0 \to A$ be a morphism to an
abelian variety.  Then the composite 
$$
\CD
G @>{\varphi_x}>> X_0 @>{f}>>  A
\endCD
$$ 
factors through a unique morphism $A(G) = G/G_{\aff} \to A$, which
must be invariant under $G_x$. Thus, $f$ factors through a unique
morphism $G/H_0 \to A$. 

Next, we claim that the Albanese morphism of $X_0$ extends to $X$. Of 
course, such an extension is unique if it exists; hence we may assume
that $X$ is quasi-projective, by Theorem
\ref{thm:quasi-projective}. Let $\psi : X \to A$ be 
a morphism as in Theorem \ref{thm:fibration}. Then $\psi$ factors
through a $G$-morphism $\alpha_{\psi}: A(X) \to A$, and the composite 
$\alpha_{\psi} \circ \alpha_{\varphi} : A(G) \to A$ is an
isogeny. Thus, $\alpha_{\varphi}$ is an isogeny, and hence the
canonical morphism $h_X : A(X)_r \to A(X)$ is an isogeny as well
(since $A(X)_r = A(X_0) = A(G)/F$). Together with Zariski's main
theorem, it follows that the rational map 
$\alpha_{X,r} : X - - \to A(X)_r$ 
is a morphism, i.e., $h_X$ is an isomorphism.

\section{Proof of Theorem \ref{thm:product}}
\label{sec:proof4}

Choose $x \in X$ so that $X \cong G/G_x$. The radical $R(G_{\aff})$
fixes some point of $X$, and is a normal subgroup of $G$; thus,
$R(G_{\aff}) \subset G_x$. By the faithfulness assumption, it follows
that $R(G_{\aff})$ is trivial, that is, $G_{\aff}$ is semi-simple. 
Moreover, the reduced connected center $C(G)^0$ satisfies 
$C(G)^0_{\aff} \subset R(G_{\aff})$. Thus, $C(G)^0$ is an abelian
variety, and hence $C(G)^0 \cap G_{\aff}$ is a finite central
subgroup scheme of $G$.

We claim that $C(G)^0 \cap G_{\aff}$ is trivial. Consider indeed the
reduced neutral component $P$ of $G_x$. Then $P \subset G_{\aff}$ by
Lemma \ref{lem:affine} (ii). Moreover, the natural map
$$
\pi: G/P \longrightarrow G/G_x \cong X
$$
is finite, and hence $G/P$ is complete. It follows that $G_{\aff}/P$
is complete as well. Thus, $P$ is a parabolic subgroup of $G_{\aff}$,
and hence contains $C(G)^0 \cap G_{\aff}$. In particular, 
$C(G)^0 \cap G_{\aff} \subset G_x$; arguing as above, this implies our
claim.

By that claim and the equality $G = C(G)^0 G_{\aff}$ (Lemma
\ref{lem:affine} (iii)), we obtain that $C(G)^0 = A(G)$ and 
$G \cong A(G) \times G_{\aff}$, where $G_{\aff}$ is semi-simple and
adjoint. Moreover,
$$
G/P \cong A(G) \times (G_{\aff}/P).
$$
Next, recall from Section \ref{sec:proof3} that the Albanese morphism
$\alpha_X$ is the natural map $G/G_x \to G/H$, where 
$H := G_x G_{\aff}$ is affine. So we may write $H = F \times G_{\aff}$,
where $F \subset A(G)$ is a finite subgroup-scheme, and hence a central
subgroup-scheme of $G$. Then 
$$
X \cong A(G) \times ^F Y
$$ 
and $A(X) \cong A(G)/F$, where 
$Y := H/G_x \cong G_{\aff}/(G_{\aff} \cap G_x)$.
We now show that $F$ is trivial; as above, it suffices to prove that
$F \subset G_x$.

Choose a maximal torus $T \subset P$ so that $x$ lies in the fixed
point subscheme $X^T$. The latter is nonsingular and stable under
$A(G)$. Moreover, the restriction
$$
\pi: (G/P)^T = A(G) \times (G_{\aff}/P)^T \longrightarrow X^T
$$
is surjective, and hence $A(G)$ acts transitively on each component of
$X^T$; the Weyl group of $(G_{\aff},T)$ permutes transitively these
components. Let $Z$ be the component containing $x$. To complete the
proof, it suffices to show that the morphism
$\alpha_X: X \to A(G)/F$ restricts to an isomorphism $Z \to A(G)/F$.

By the Bialynicki-Birula decomposition (see \cite{Bi73}), $Z$ admits a
$T$-stable open neighborhood $U \subset X$ together with a
$T$-equivariant retraction $\rho : U \to Z$, a locally trivial
fibration in affine spaces. It follows that the Albanese morphisms of
$U$ and $Z$ satisfy $\alpha_U = \alpha_Z \circ \rho$. Since $\alpha_U$
is the restriction of $\alpha_X$, we obtain that 
$\alpha_Z = \alpha_X \vert_Z$. On the other hand, $\alpha_Z$ is the
identity, since $Z$ is homogeneous under $A(G)$; this completes the
argument.

\section{Examples}
\label{sec:examples}

\begin{example}\label{ex:several}
Let $A$ be an abelian variety. Choose distinct finite subgroups
$F_1, \ldots, F_n \subset A$ such that 
$F_1 \cap \cdots \cap F_n = \{ 0 \}$ (as schemes), and let
$$
X := A/F_1 \times \cdots \times A/F_n.
$$
Then $A$ acts faithfully on $X$ via simultaneous translation on all
factors, and each projection $p_i : X \to A/F_i$ satisfies the
assertion of Theorem \ref{thm:fibration}. Since the fibers of $p_i$
are varieties, this morphism admit no non-trivial factorization
through a $A$-morphism $\psi: X \to A/F$, where $F$ is a finite
subgroup scheme of $A$.

Assume that the only endomorphisms of $A$ are the multiplications by
integers, and the orders of $F_1, \ldots,F_n$ have a non-trivial
common divisor. Then there exists no $A$-morphism $\psi : X \to A$
(otherwise, we would obtain endomorphisms 
$u_1, \ldots, u_n$ of $A$ such that each $u_i$ is zero on $F_i$ and 
$u_1 + \cdots + u_n = \id_A$, which contradicts our
assumptions). Equivalently, $A$ is not a direct summand of $X$ for its 
embedding via the diagonal map $a \mapsto (a + F_1,\ldots,a + F_n)$.
Thus, there exists no universal morphism $\psi$ satisfying the
assertions of Theorem \ref{thm:fibration}.

\end{example}

\begin{example}\label{ex:raynaud1}
Following \cite[XIII 3.1]{Ra70}, we construct a complete equivariant
embedding $X$ of a commutative non-affine group $G$, for which the
assertions of Theorems 1--3 do not hold. (Of course, $X$ will be
non-normal).

Consider an abelian variety $A$, and a point $a \in A$.
Let $X_a$ be the scheme obtained from $A \times \bP^1$
by identifying every point $(x,0)$ with $(x + a, \infty)$.
(That $X_a$ is indeed a scheme follows e.g. from the main result of
\cite{Fe03}.) Then $X_a$ is a complete variety, and the canonical map
$$
f_a: A \times \bP^1 \longrightarrow X_a
$$ 
is the normalization. In particular, $X_a$ is \emph{weakly normal},
i.e., every finite birational bijective morphism from a variety to
$X_a$ is an isomorphism.

The connected algebraic group 
$$
G := A \times \bG_m
$$ 
acts faithfully on $A \times \bP^1$ via 
$(x,t) \cdot (y,u) = (x + y, tu)$, and this induces a faithful
$G$-action on $X_a$ such that $f_a$ is equivariant. Moreover, $X_a$
consists of two $G$-orbits: the open orbit,
$$
f_a(A \times (\bP^1 \setminus \{ 0,\infty \})) \cong G,
$$
and the closed orbit, 
$$
f_a(A \times \{ 0\})  = f_a(A \times \{ \infty \} ) 
\cong A \cong G/\bG_m,
$$
which is the singular locus of $X_a$. 

The projection $A \times \bP^1 \to A$ induces a morphism 
$$
\alpha_a : X_a \longrightarrow A/a,
$$ 
where $A/a$ denotes the quotient of $A$ by the closure of the subgroup
generated by $a$. We claim that $\alpha_a$ is the Albanese morphism:
indeed, any morphism $\beta : X_a \to B$, where $B$ is an abelian
variety, yields a morphism $\beta \circ f_a : A \times \bP^1 \to B$. By
rigidity (see e.g. \cite[Cor.~2.5, Cor.~3.9]{Mi86}), there exists a
morphism $\gamma : A \to B$ such that 
$$
(\beta \circ f_a)(y,z) = \gamma(y)
$$ 
for all $y \in A$ and $z \in \bP^1$. Conversely, a morphism 
$\gamma : A \to B$ yields a morphism $\beta : X_a \to B$ if and only
if $\gamma(y + a) = \gamma(y)$ for all $y \in A$, which proves our claim.

Also, note that $\alpha_G : G \to A(G)$ is just the projection 
$A \times \bG_m \to A$. Thus, the kernel of the homomorphism 
$A(G) \to A(X_a)$ is the closed subgroup generated by $a$.

In particular, if the order of $a$ is infinite, then $X_a$ does not
admit any $G$-morphism to a finite quotient of $A$, so that the
assertions of Theorems \ref{thm:fibration} and
\ref{thm:albanese} are not satisfied. Moreover, $X_a$ is not
projective (as follows e.g. from Lemma \ref{lem:equivalent} applied to
the action of $A$), so that the assertion of Theorem
\ref{thm:quasi-projective} does not hold as well. 

On the other hand, if $a$ has finite order $n$, then the fibers of
$\alpha_a$ are unions of $n$ copies of $\bP^1$ on which the origins
are identified. In that case, $X_a$ is projective, but does not
satisfy the assertions of Theorem \ref{thm:albanese} when $n\ge 2$. 

We may also consider $X_a$ as an $A$-variety via the inclusion 
$A \subset G$. The projection $A \times \bP^1 \to \bP^1$ induces an
$A$-invariant morphism
$$
\pi: X_a \longrightarrow C,
$$ 
where $C$ denotes the nodal curve obtained from $\bP^1$ by identifying
$0$ with $\infty$; one checks that $\pi$ is an $A$-torsor. If $a$ has
infinite order, then $X_a$ is not covered by $A$-stable quasi-projective 
open subsets: otherwise, by Lemma \ref{lem:equivalent} again, any such
subset $U$ would admit an $A$-morphism to a finite quotient of $A$,
i.e., the kernel of the $A$-action on $A(U)$ would be finite. However,
when $U$ meets the singular locus of $X_a$, one may show as above that
the Albanese morphism of $U$ is the restriction of $\alpha_a$, which
yields a contradiction.
\end{example}

\begin{example}\label{ex:raynaud2}
Given any abelian variety $A$, we construct (after 
\cite[XII 1.2]{Ra70}) an $A$-torsor $\pi: X \to Y$ such that $A$ acts
trivially on $A(X)$.

Let $X$ be the scheme obtained from 
$A \times A \times \bP^1$ by identifying every point 
$(x,y,0)$ with $(x+y, y, \infty)$. Then again, $X$ is a complete
variety, and the canonical map 
$$
f: A \times A \times \bP^1 \longrightarrow X
$$
is the normalization; $X$ is weakly normal but not normal.

Let $A$ act on $A \times A \times \bP^1$ via translation on the first
factor; this induces an $A$-action on $X$ such that $f$ is
equivariant. Moreover, the projection
$p_{23}: A \times A \times \bP^1 \to A \times \bP^1$
induces an $A$-invariant morphism 
$$
\pi: X \longrightarrow A \times C,
$$
where $C$ denotes the rational curve with one node, as in the above
example. One checks that $\pi$ is an $A$-torsor, and the Albanese
morphism $\alpha_X$ is the composite of $\pi$ with the projection 
$A \times C \to A$. Thus, $\alpha_X$ is $A$-invariant.

So Theorem \ref{thm:fibration} does not hold for the $A$-variety $X$. 
Given any $A$-stable open subset $U \subset X$ which meets the
singular locus, one may check as above that $\alpha_U$ is
$A$-invariant; as a consequence, Theorem \ref{thm:quasi-projective}
does not hold as well.
\end{example}

\begin{example}\label{ex:raynaud3}
Given an elliptic curve $E$, we construct an example of an $E$-torsor
$\pi: X \to Y$, where $Y$ is a normal affine surface, and the Albanese
variety of $X$ is trivial. (In particular, $X$ is a normal $E$-variety
which does not satisfy the assertions of Theorem \ref{thm:fibration}).
For this, we adapt a construction in \cite[XIII 3.2]{Ra70}.

Let $\dot E := E \setminus \{ 0 \}$ and 
$\dot \bA^1 := \bA^1 \setminus \{ 0 \}$. We claim that there exist a
normal affine surface $Y$ having exactly two singular points
$y_1,y_2$, and a morphism
$$
f :  \dot E \times \dot \bA^1 \longrightarrow Y
$$
such that $f$ contracts $\dot E \times \{ 1 \}$ to $y_1$,
$\dot E \times \{ -1 \}$ to $y_2$, 
and restricts to an isomorphism on the open subset
$\dot E \times (\dot \bA^1 \setminus \{1,-1\})$.

Indeed, we may embed $E$ in $\bP^2$ as a cubic curve with homogeneous
equation $F(x,y,z) = 0$, such that the line $(z=0)$ meets $C$ with
order $3$ at the origin. Then $\dot E$ is identified with
the curve in $\bA^2$ with equation $F(x,y,1) = 0$. Now one readily
checks that the claim holds for the surface
$$
Y \subset \bA^2 \times \dot \bA^1
$$
with equation $F(x,y,z^2 - 1) = 0$, and the morphism
$$
f : (x,y,z) \longmapsto (x(z^2 - 1), y(z^2 - 1), z).
$$
The singular points of $Y$ are $y_1 = (0,0,1)$ and $y_2 = (0,0,-1)$.

Next, let $U_1 = Y \setminus \{ y_2 \}$, 
$U_2 := Y \setminus \{ y_1 \}$, and $U_{12} := U_1 \cap U_2$. 
Then $U_{12}$ is nonsingular and contains an open subset 
isomorphic to $\dot E \times (\dot \bA^1 \setminus \{ 1, -1 \})$.
Thus, the projection 
$\dot E \times (\dot \bA^1 \setminus \{ 1, -1 \}) \to \dot E$
extends to a morphism
$$ 
p : U_{12}  \longrightarrow E.
$$
We may glue $E \times U_1$ and $E \times U_2$ along $E \times U_{12}$
via the automorphism 
$$
(x,y) \longmapsto (x + p(y),y)
$$
to obtain a torsor $\pi: X \to Y$ under $E$. Arguing as in Example
\ref{ex:raynaud1}, one checks that the Albanese variety of $X$ is
trivial.
\end{example}

\begin{example}\label{ex:igusa}
Following Igusa (see \cite{Ig55}), we construct nonsingular projective 
varieties for which the tangent bundle is trivial, and the Albanese
morphism is a non-trivial fibration.

We assume that the ground field $k$ has characteristic $2$, and
consider two abelian varieties $A$, $B$ such that $A$ has a point $a$
of order $2$. Let $X$ be the quotient of $A \times B$ by the involution 
$$
\sigma : (x,y) \longmapsto (x + a, -y)
$$
and denote by $\pi : A \times B \to X$ the quotient morphism. Since
$\sigma$ has no fixed point, $X$ is a nonsingular projective variety,
and the natural map between tangent sheaves
$\cT_{A \times B} \to \pi^*(\cT_X)$ is an isomorphism. 
This identifies $\Gamma(X, \cT_X)$ with the invariant subspace
$$
\Gamma(A \times B, \cT_{A \times B})^{\sigma} = \Lie(A \times B)^{\sigma}.
$$
Moreover, the action of $\sigma$ on $\Lie(A \times B)$ is trivial, by
the characteristic $2$ assumption. This yields
$$
\Gamma(X, \cT_X) \cong \Lie(A \times B).
$$
As a consequence, the natural map 
$\cO_X \otimes \Gamma(X,\cT_X) \to \cT_X$ is an isomorphism, since this
holds after pull-back via $\pi$. In other words, the tangent sheaf of
$X$ is trivial.

The action of $A$ on $A \times B$ via translation on the first factor
induces an $A$-action on $X$, which is easily seen to be faithful. 
Moreover, the first projection $A \times B \to A$ induces an
$A$-morphism 
$$
\psi: X \longrightarrow A/a,
$$
a homogeneous fibration with fiber $B$. One checks as in Example
\ref{ex:raynaud1} that $\psi$ is the Albanese morphism of $X$.

Let $G$ denote the reduced neutral component of the automorphism group 
scheme $\Aut(X)$. We claim that the natural map $A \to G$ is an
isomorphism. Indeed, $G$ is a connected algebraic group, equipped with
a morphism $G \to A/a$ having an affine kernel (by the
Nishi--Matsumura theorem) and such that the composite 
$A \to G \to A/a$ is the quotient map. Thus, $G = A G_{\aff}$, and
$G_{\aff}$ acts faithfully on the fibers of $\psi$, i.e., on $B$. This
implies our claim.

In particular, $X$ is not homogeneous, and the fibration $\psi$ is
non-trivial (otherwise, $X$ would be an abelian variety). In fact, the
sections of $\psi$ correspond to the $2$-torsion points of $B$, and
coincide with the local sections; in particular, $\psi$ is not locally
trivial.
\end{example}

\begin{example}\label{ex:split}
We assume that $k$ has characteristic $p > 0$. Let $X$ be the
hypersurface in $\bP^n \times \bP^n$ with bi-homogeneous equation
$$
f(x_0, \ldots, x_n, y_0, \ldots, y_n) := \sum_{i=0}^n x_i^p \, y_i \, =0,
$$
where $n \ge 2$. The simple algebraic group $G := \PGL(n+1)$ 
acts on $\bP^n \times \bP^n$ via
$$
[A] \cdot ([x,],[y]) = ([A x], [F(A^{-1})^T \, y]),
$$
where $A \in \GL(n+1)$, and $F$ denotes the Frobenius endomorphism of
$\GL(n+1)$ obtained by raising matrix coefficients to their $p$-th
power. This induces a $G$-action on $X$ which is faithful and
transitive, but not smooth; the isotropy subgroup scheme of any point
$\xi = ([x],[y])$ is the intersection of the parabolic subgroup
$G_{[x]}$ with the non-reduced parabolic subgroup scheme $G_{[y]}$. So
$X$ is a variety of unseparated flags.

Denote by 
$$
\pi_1, \pi_2 : X \longrightarrow \bP^n
$$
the two projections. Then $X$ is identified via $\pi_1$ to the 
projective bundle associated to $F^*(\cT_{\bP^n})$ (the pull-back of
the tangent sheaf of $\bP^n$ under the Frobenius morphism). In
particular, $\pi_1$ is smooth. Also, $\pi_2$ is a homogeneous
fibration, and $(\pi_2)_*\cO_X = \cO_{\bP^n}$, but $\pi_2$
is not smooth.

Let $\cT_{\pi_1}$, $\cT_{\pi_2}$ denote the relative tangent sheaves
(i.e., $\cT_{\pi_i}$ is the kernel of the differential 
$d\pi_i : \cT_X \to \pi_i^*(\cT_{\bP^n})$). We claim that 
\begin{equation}\label{eqn:split}
\cT_X = \cT_{\pi_1} \oplus \cT_{\pi_2} \,;
\end{equation}
in particular, $\cT_{\pi_2}$ has rank $n-1$. Moreover, $\cT_{\pi_2}$
is the subsheaf of $\cT_X$ generated by the global sections. 

For this, consider the natural map 
$$
\op_X: \Lie(G) \longrightarrow \Gamma(X,\cT_X)
$$
and the induced map of sheaves
$$
op_X: \cO_X \otimes \Lie(G) \longrightarrow \cT_X.
$$
For any $\xi = ([x], [y]) \in X$, the kernel of the map 
$op_{X,\xi} : \Lie(G) \to \cT_{X,\xi}$ 
is the isotropy Lie algebra $\Lie(G)_{\xi}$, that is, $\Lie(G)_{[x]}$
(since $\Lie(G)$ acts trivially on the second copy of $\bP^n$). Thus, 
$(d\pi_1)_{\xi}: \cT_{X,\xi} \longrightarrow \cT_{\bP^n,[x]}$
restricts to an isomorphism $\Ima(op_{X,\xi}) \cong \cT_{\bP^n,[x]}$.
In other words, $d\pi_1$ restricts to an isomorphism
$\Ima(op_X) \cong \pi_1^*(\cT_{\bP^n})$.
So $\pi_1$ is the Tits morphism of $X$. Moreover, we have a decomposition
$$
\cT_X = \cT_{\pi_1} \oplus \Ima(op_X),
$$
and $\Ima(op_X) \subset \cT_{\pi_2}$. Since 
$\cT_{\pi_1} \cap \cT_{\pi_2} = 0$, this implies (\ref{eqn:split}) and
the equality $\Ima(op_X) = \cT_{\pi_2}$. 

To complete the proof of the claim, we show that $\op_X$ is an
isomorphism. Consider indeed the bi-homogeneous coordinate ring
$$
k[x_0, \ldots, x_n,y_0,\ldots,y_n]/(x_0^p \, y_0 + \cdots + x_n^p \, y_n)
$$
of $X$. Its homogeneous derivations of bi-degree $(0,0)$ are those
given by $x \mapsto Ax$, $y \mapsto ty$, where $A$ is an 
$n+1 \times n+1$ matrix, and $t$ a scalar; this is equivalent to our
assertion.

Also, $X$ admits no non-trivial decomposition into a direct product,
by \cite[Cor.~6.3]{Sa03}. This may be seen directly: 
if $X \cong X_1 \times X_2$, then the $G$-action on $X$ induces
actions on $X_1$ and $X_2$ such that both projections are equivariant
(as follows e.g. by linearizing the pull-backs to $X$ of ample
invertible sheaves on the nonsingular projective varieties $X_1$,
$X_2$). So $X_1 \cong G/H_1$ and $X_2 \cong G/H_2$, where $H_1$, $H_2$
are parabolic subgroup schemes of $G$. Since the diagonal $G$-action
on $G/H_1 \times G/H_2$ is transitive, it follows that $G = P_1 P_2$,
where $P_i$ denotes the reduced scheme associated to $H_i$ (so that
each $P_i$ is a proper parabolic subgroup of $G$). But the simple
group $G$ cannot be the product of two proper parabolic subgroups, a
contradiction.

Since $X$ is rational and hence simply-connected, this shows that 
a conjecture of Beauville (relating decompositions of the tangent
bundle of a compact K\"ahler manifold to decompositions of its
universal cover, see \cite{Be00}) does not extend to nonsingular
projective varieties in positive characteristics.

More generally, let $G$ be a simple algebraic group in characteristic
$p \ge 5$, and $X \cong G/G_x$ a complete variety which is homogeneous
under a faithful $G$-action. By \cite{We93}, there exists a unique
decomposition
$$
G_x = \bigcap_{i=1}^r P_i \, G_{n_i} \, ,
$$
where $P_1,\ldots,P_r$ are pairwise distinct maximal parabolic
subgroups of $G$, and $(n_1,  \ldots,  n_r)$ is an increasing sequence
of non-negative integers. Here $G_n$ denotes the $n$-th Frobenius
kernel of $G$. Since $G$ acts faithfully on $X$, we must have 
$n_1 = 0$. Let 
$$
Q_1 := \bigcap_{i, n_i = 0} P_i \, , \quad 
Q_2 := \bigcap_{i, n_i \ge 1} P_i \, G_{n_i}\, . 
$$
Then $Q_1$ is a parabolic subgroup of $G$, $Q_2$ is a parabolic
subgroup scheme, $G_x = Q_1 \cap Q_2$, and $\Lie(G_x) = \Lie(Q_1)$. 
Thus, the Tits morphism of $X$ is the canonical map
$$
\pi_1 : G/G_x \to G/Q_1,
$$ 
and $d\pi_1$ restricts to an isomorphism 
$\Ima(op_X) \cong \pi_1^*(\cT_{G/Q_1})$. This implies the
decomposition $\cT_X = \Ima(op_X) \oplus \cT_{\pi_1}$, which is
non-trivial unless $G_x = Q_1$ (i.e., $G_x$ is reduced). Moreover, 
$\Ima(op_X) \subset \cT_{\pi_2}$ (where $\pi_2 : G/G_x \to G/Q_2$
denotes the canonical map), since $G_1$ acts trivially on $G/Q_2$. 
As $\pi_1 \times \pi_2$ is a closed immersion, it follows again that
$$
\cT_X = \cT_{\pi_1} \oplus \cT_{\pi_2} \quad \text{and} \quad  
\Ima(op_X) = \cT_{\pi_2}.
$$
On the other hand, the variety $X$ is indecomposable, since $G$ is
simple (see \cite[Cor.~6.3]{Sa03} again, or argue directly as
above). This yields many counter-examples to the analogue of the
Beauville conjecture. 
\end{example}


\begin{thebibliography}{100} 

\bibitem[Ak]{Ak95} D.~Akhiezer,
\emph{Lie group actions in complex analysis},
Aspects of Mathematics {\bf E 27}, Vieweg, Braunschweig/Wiesbaden, 1995.

\bibitem[Al]{Al02} V.~Alexeev,
\emph{Complete moduli in the presence of semi-abelian group action},
Ann. of Math. {\bf 155} (2002), 611--718.

\bibitem[Be]{Be00} A.~Beauville, 
\emph{Complex manifolds with split tangent bundle}, 
in: Complex analysis and algebraic geometry, 61--70, de Gruyter,
Berlin, 2000.

\bibitem[Bi]{Bi73} A.~Bialynicki-Birula, 
\emph{Some theorems on actions of algebraic groups},
Ann. Math. (2) {\bf 98} (1973), 480--497.

\bibitem[BLR]{BLR90} S.~Bosch, W.~L\"utkebohmert and M.~Raynaud,
\emph{N\'eron Models}, Ergeb. der Math. {\bf 21}, Springer-Verlag,
Berlin, 1990.

\bibitem[Br]{Br06} M.~Brion,
\emph{Log homogeneous varieties},
arXiv: math.AG/069669, to appear in the proceedings of the VI
Coloquio Latinoamericano de \'Algebra (Colonia, Uruguay, 2005).

\bibitem[BR]{BR06} M.~Brion and A.~Rittatore,
\emph{The structure of normal algebraic monoids},
arXiv: math.AG/0610351, to appear in Semigroup Forum. 

\bibitem[Co]{Co02} B.~Conrad,
\emph{A modern proof of Chevalley's theorem on algebraic groups},
J. Ramanujam Math. Soc. {\bf 17} (2002), 1--18.

\bibitem[DG]{DG70} M.~Demazure and P.~Gabriel,
\emph{Groupes alg\'ebriques}, North Holland, Amsterdam, 1970.

\bibitem[Fe]{Fe03} D.~Ferrand,
\emph{Conducteur, descente et pincement},
Bull. Soc. math. France {\bf 131} (2003), 553--585.

\bibitem[HL]{HL93} W.~Haboush and N.~Lauritzen, 
\emph{Varieties of unseparated flags,}
in: Linear algebraic groups and their representations (Los Angeles,
CA, 1992), 35--57, Contemp. Math. {\bf 153}, Amer. Math. Soc.,
Providence, RI, 1993. 

\bibitem[Ha]{Ha77} R.~Hartshorne,
\emph{Algebraic Geometry}, Grad. Text Math. {\bf 52}, Springer--Verlag,
New York, 1977.

\bibitem[Ig]{Ig55} J.-i.~Igusa,
\emph{On some problems in abstract algebraic geometry},
Proc. Nat. Acad. Sci. U. S. A. {\bf 41} (1955), 964--967.

\bibitem[KSS]{KSS89} H.~Kraft, P.~Slodowy and T.~A.~Springer (eds.),
\emph{Algebraic Transformation Groups and Invariant Theory},
DMV Seminar Band {\bf 13}, Birkh\"auser, Basel, 1989.

\bibitem[Ma]{Ma63} H.~Matsumura,
\emph{On algebraic groups of birational transformations},
Atti Accad.~Naz.~Lincei Rend.~Cl.~Sci.~Fis.~Mat.~Natur. (8)
{\bf 34} (1963), 151--155.

\bibitem[Mi]{Mi86} J.~S.~Milne,
\emph{Abelian Varieties}, in:
Arithmetic Geometry (G.~Cornell and J.~H.~Silverman, eds.), 103--150,
Springer--Verlag, New York, 1986. 

\bibitem[Muk]{Mu78} S.~Mukai,
\emph{Semi-homogeneous vector bundles on Abelian varieties},
J. Math. Kyoto Univ. {\bf 18} (1978), 239--272.

\bibitem[MFK]{MFK94} D.~Mumford, J.~Fogarty and F.~Kirwan,
\emph{Geometric Invariant Theory}, 3rd edition, 
Springer-Verlag, New York, 1994.

\bibitem[Mum]{Mu70} D.~Mumford,
\emph{Abelian Varieties},
Oxford University Press, Oxford, 1970.

\bibitem[PV]{PV94} V.~L.~Popov and E.~B.~Vinberg,
\emph{Invariant Theory}, in: Algebraic Geometry IV,
Encycl. Math. Sci. {\bf 55}, Springer-Verlag, New York, 1994.

\bibitem[Ra]{Ra70} M.~Raynaud,
\emph{Faisceaux amples sur les sch\'emas en groupes et les espaces
homog\`enes}, Lecture Notes in Math. {\bf 119}, Springer-Verlag, New
York, 1970.

\bibitem[Sa]{Sa03} C.~Sancho de Salas,
\emph{Complete homogeneous varieties: structure and classification},
Trans. Amer. Math. Soc. {\bf 355} (2003), 3651--3667.

\bibitem[Se1]{Se58a} J.-P.~Serre, 
\emph{Espaces fibr\'es alg\'ebriques,}
S\'eminaire C.~Chevalley (1958), Expos\'e No.~1, 
Documents Math\'ematiques {\bf 1}, Soc. Math. France, Paris, 2001.

\bibitem[Se2]{Se58b} J.-P.~Serre, 
\emph{Morphismes universels et vari\'et\'e d'Albanese},
S\'eminaire C.~Chevalley (1958--1959), Expos\'e No.~10, 
Documents Math\'ematiques {\bf 1}, Soc. Math. France, Paris, 2001.

\bibitem[Su]{Su74} H.~Sumihiro,
\emph{Equivariant completion}, 
J. Math. Kyoto Univ. {\bf 14} (1974), 1--28.

\bibitem[We]{We93} C.~Wenzel, 
\emph{Classification of all parabolic subgroup-schemes of a reductive
linear algebraic group over an algebraically closed field},
Trans. Amer. Math. Soc. {\bf 337} (1993), 211--218. 

\end{thebibliography}
\end{document}